\documentclass[preprint,12pt]{elsarticle}
\usepackage{graphicx}

\usepackage{setspace}
\usepackage[margin=1.5in]{geometry}    
\usepackage[utf8]{inputenc}
\usepackage{algorithm}
\usepackage[noend]{algpseudocode}

\makeatletter
\def\BState{\State\hskip-\ALG@thistlm}
\makeatother
\usepackage{enumitem} 
\usepackage{tablefootnote}
\usepackage{caption}
\usepackage{graphicx}
\usepackage{booktabs}
\usepackage{amsfonts}
\usepackage{amsthm}
\usepackage{amsmath,bm}
\usepackage{amssymb}
\usepackage{mathtools}
\usepackage{subfig}
\usepackage{multirow}
\usepackage{bbm}
\usepackage{hyperref}

\usepackage{array}

\usepackage[colorinlistoftodos]{todonotes}
\usepackage{comment}

\usepackage{mathrsfs}

\usepackage{epstopdf}
\usepackage{enumitem}  

\journal{Transportation Research Part C}

\begin{document}

\begin{frontmatter}


\title{A Two-Stage Stochastic Programming Model for Car-Sharing Problem using Kernel Density Estimation}

\author[1,3]{Xiaoming Li}
\author[1]{Chun Wang}
\author[2]{Xiao Huang}

\address[1]{<Department of Information and System Engineering, Concordia University, Montreal, QC, Canada>}
\address[2]{<John Molson School of Business, Concordia University, Montreal, QC, Canada>}
\address[3]{<School of Computer, Shenyang Aerospace University, Shenyang, Liaoning, China>}

\begin{abstract}
Car-sharing problem is a popular research field in sharing economy. In this paper, we investigate the car-sharing re-balancing problem under uncertain demands. An innovative framework that integrates a non-parametric approach - kernel density estimation (KDE) and a two-stage stochastic programming (SP) model is proposed. Specifically, the probability distributions are derived from New York taxi trip data sets by KDE, which are used as the input uncertain parameters for SP. Additionally, the car-sharing problem is formulated as a two-stage SP model which aims to maximize the overall profit. Meanwhile, a Monte Carlo method called sample average approximation (SAA) and Benders decomposition algorithm are introduced to solve the large-scale optimization model. Finally, the experimental validations show that the proposed framework outperforms the existing works in terms of outcomes.
\end{abstract}


\end{frontmatter}

\setlength{\baselineskip}{20pt}

\section{Introduction}
\label{S:1}

\par
Car-Sharing which was coined in the middle of 20th century\cite{shaheen1998carsharing} is increasing sharply in many cities now. This trend becomes much more popular since people may benefit a lot from the the sharing system, such as saving parking lots, reducing the traffic congestion and air pollution\cite{bruglieri2018two}. To use the car-sharing service, normally, a customer can reserve the vehicle by phone or Internet. Once approved, the reserved vehicle is assigned to the customer who picks it up at an appointed time and leaves it at a specific car-sharing location, which may be the same as the pick-up point (one-way car-sharing systems) or anywhere in a specified zone (free-floating car-sharing systems) \cite{vosooghi2017critical, illgen2018literature}. Typically, car-sharing systems are financed by public and /or private entities and managed by service provides, who are involved in strategic, tactical, and operational decision-making. Strategic decisions can be include determining the number, location, and capacity of stations for car rental and return, whereas tactical decisions can include allocation decisions. Daily, operational decisions include determining how to periodically re-distribute cars to station. 

\par
Nowadays, With the rapid development of transportation in cities, a huge amount of data is generated every day. Thanks to the emerging technologies such as wireless sensor network (WSN), cloud computing and Big Data, which make it possible to collect, store and analyze the data in an effective and efficient way. However, increasing data brings new challenges to traditional car-sharing optimization issues.

\par
In order to solve the problem, several optimization approaches are under investigation including complicating determinate modeling and optimization under uncertainty. 

\par
For deterministic model, Gambella et al. \cite{gambella2018optimizing} propose an MIP model along with two heuristic algorithms to optimize electric vehicle relocation problem, Huang et al.\cite{huang2018solving} investigate one-way station-based relocation considering non-linear demand, an mixed integer non-linear model is proposed, Xu et al.\cite{xu2018electric} study the electric vehicle fleet size and trip pricing (EVFS\&TP) problem for one-way car-sharing services by taking into account the necessary practical requirements of vehicle relocation and personnel assignment. A mixed-integer nonlinear and nonconvex programming model is developed. Li et al. \cite{li2016adaptive} focus on the Share-a-Ride Problem (SARP) aiming at maximizing the profit of serving a set of passengers and parcels using a set of homogeneous vehicles. An  adaptive large neighborhood search heuristic algorithm is devised. Zhao et al. \cite{zhao2018integrated} devise an integrated framework to minimize the total cost, including the EV and staff investment, EV re-balancing and staff relocation costs. The model is reformulated and solved by Lagrangian relaxation approach. Boyacı et al. \cite{boyaci2015optimization} explore one-way vehicle-sharing systems that is taking vehicle relocation and electric vehicle charging requirements into consideration. A multi-objective optimization model is developed and solved by branch-and-bound. 

\par
For stochastic programming model, Brandstätter et al.\cite{brandstatter2017determining} solve strategic optimization problems of car-sharing systems that utilize electric cars by a two-stage stochastic programming model. Also, the heuristic algorithm is used to tackle large-scale instances. Biondi et al. \cite{biondi2016optimal} explore to optimize car-sharing system with uncertain demands from the perspective of queue theory. Fan et al. \cite{fan2004multi} consider the stochastic dynamic vehicle allocation problem (SDVAP), a multi-stage stochastic programming model is formulated to maximize profits and to manage fleets of vehicles in both time and space. Later, they develop a stochastic programming model to optimize strategic allocation of vehicles for one-way car-sharing systems under demand uncertainty
\cite{fan2014optimizing}. Cavagnini et al. \cite{cavagnini2018two} propose a bike-sharing system which composes one depot and multiple capacitated stations.

\par
Although the aforementioned work handle with car-sharing problem from different perspectives, among these works, most of the works consider modeling car-sharing problem in the deterministic way without involving any uncertain parameters such as demand, supply, travelling time. Only a few of them use optimization under certainty techniques to solve car-sharing problem, none of them utilize the accurate probability information from historical data. Meanwhile, most mathematical models that are formulated based on SP are assumed that the probability distribution is known with a specific type. However, in the real historical data, the probability distribution information may contain many even infinite parameters which cannot be described by simple known distribution such as Gaussian distribution. We will further discuss this topic in the numerical experiment section.

\par
There are several ways to hedge against uncertainty using optimization techniques. In stochastic programming(SP)\cite{birge2011introduction}, uncertainty is modeled through discrete or continuous probability functions, in other words, SP models heavily rely on probability information from historical data. In fuzzy programming (FP)\cite{zimmermann1978fuzzy}, uncertainty parameters are considered as fuzzy numbers and constraints are treated as fuzzy sets. In robust optimization(RO)\cite{ben2009robust}, uncertainty is described in a particular set called uncertain set. In distibtuionally robust optimization(DRO)\cite{delage2010distributionally}, uncertainty is formulated by an ambiguity set which includes a family of probability distributions. In our scope, we are primarily interested in extracting exact probability distribution information from historical data.  To this end, we come to consider using two-stage SP to solve the car-sharing problem. 

\par
In order to overcome the issue aforementioned, we consider to utilize related machine learning approach to make the SP model more practical. Recently, integrating machine learning (ML) with optimization techniques becomes the trend in operational research (OR) community\cite{bengio2018machine}, \cite{larsen2018predicting}. A few researchers attempted to leverage the advantages of ML make optimization more realistic, especially, when it is applied in big data and data-driven optimization\cite{ning2018data},\cite{shang2018distributionally}. In our work, we will follow the trend to solve car-sharing problem. Specifically, we proposed a framework that involves two major components. In ML part, we utilize the non-parametric approach - kernel density estimation to extract more accurate probability distribution from historical data, while in OR part, stochastic programming models are constructed based on those parameters. To our best knowledge, the proposed framework is the first one to solve car-sharing problem under demand uncertainty. The contribution of this work can be summarized as follows.
\begin{enumerate}[label=(\roman*)]
    \item We consider using the non-parametric approach kernel density estimation to extract the arbitrary probability distribution of user demands from historical data on New York taxi trip data set.
    \item A two-stage stochastic programming model using the aforementioned probability distribution information is proposed to formulate car-sharing problem.
    \item Integrating sample average approximation method with Benders decomposition algorithm to solve the two-stage stochastic programming model.
\end{enumerate}

\par
The rest of the paper is organized as follows. The problem description is discussed in section 2. Section 3 investigates some related literature, the methodology is explored in section 4. In section 5, both deterministic and two-stage stochastic models are designed. While section 6 describes the framework which involves sample average approximation (SAA) method and Benders decomposition algorithm. Data prepossession and numerical experiment are presented in section 7. Finally, we conclude our work and propose future work in section 8. 
\section{Problem Description}
\label{S:2}

In this article, we address the car-sharing problem with the demands under uncertainty using a two-stage stochastic programming model. The objective is to make the maximize overall profit, which involves total revenue, holding costs at each location and moving costs between locations. Generally, we study a car-sharing system managed by a service provider wherein the decision-making is centralized. The problem can be stated as follows: the decision making can be divided into two stages. During the first stage, at the beginning of the day (e.g. at 0 AM), the number of vehicles at each location must be determined. During the second stage, after the real demand revealed (e.g. no new orders for today accepted), the truck carriers must decide how many vehicles to relocate between locations. 
\par
In this work, the most critical concern for car-sharing problem is the way of modeling uncertainty. For convenience, only customer demands considered as uncertainty parameters. Normally, in SP paradigm uncertain parameters are modeled as random variables with specific probability distributions which are extracted from historical data. Unlike the existing works which assume the uncertainty parameters conform a known probability distribution such as Gaussian, Poisson, log-normal distributions, in our work, the uncertainty parameters follow any types of distributions or non-parametric distributions.

\par
In the two-stage SP model, all the decision variable are divided into two groups: the first stage decision variables (or here-and-now) which should be determined before the real demands revealed, and the second stage decision variables (or wait-and-see) which are determined after the real demands realized. Based on the problem statement, a group of assumption are made as follows.
\begin{enumerate}
    \item Allocation resource is limited, cannot satisfy all the demands,
    \item Holding costs incur depends on the number of vehicles and locations,
    \item Moving cost depends on the specific route,
    \item The car-sharing service must finish within one day.
\end{enumerate}
\section{Model Formulation}
In this section, we will discuss car-sharing model formulations include deterministic model and two-stage SP counterpart. It is worth noting that probability distributions are required for SP model. Unlike most existing works which assume that the probability distribution of uncertain parameters are known, the probability distribution information in our work is obtained from non-parametric learning approach - kernel density estimation. For clarity, the notations are listed in below. 
\\
\underline{\textbf{Indices/Sets}} \\
$i, j \in \mathscr{R}$ regional origins and/or destinations 
\\
$s \in \mathscr{S}$ The set of scenarios     
\\
\underline{\textbf{Parameters}} \\
$h_{i}$ = holding cost at location $i$. \\
$t_{i,j}$ = moving cost from location $i$ to location $j$. \\
$d_{}^{avg}$ = the average demand of location $j$. 
\\
\underline{\textbf{Decision Variables}} \\
$x_{i}$ = first-stage decision variable which denotes the number of vehicles at location $i$. 
\\
\underline{\textbf{Random Variables (for stochastic programming model)}} 
\\
$\widetilde{d}_{i}$ = random demands which denotes the number of cars that will be picked up by customers at location $i$.
\\
$p^{s}$ = the probability of scenario s.
\\
$y_{i,j}^{s}$ = the second-stage decision variable which denotes the number of vehicles moving from location i to location j under scenario s.

\subsection{Deterministic Model}
In the deterministic model, we consider to allocate the limited vehicles to different locations in order to maximize the overall profit. For convenience, we consider using the average demands. The deterministic model for car-sharing problem can be formulated as follows.

\begin{equation}
\label{one_1}
\text{max} \left[\sum_{i \in \mathscr{R}} \min \left(x_{i}+\sum_{j \in \mathscr{R}} y_{i,j}, d_{i}^{avg}\right) * r_{i} - \sum_{i \in \mathscr{R}}\left(h_{i} * x_{i}+\sum_{j \in \mathscr{R}} t_{i,j} * y_{i,j} \right) \right]
\end{equation}

\textbf{s.t.}
\begin{equation}
\label{one_2}
\sum_{i \in \mathscr{R}} x_{i} \leqslant \mathcal{C},
\end{equation}

\begin{equation}
\label{one_3}
\sum_{j \in \mathscr{R}} y_{i,j} \leq x_{i} \quad \forall i \in \mathscr{R},
\end{equation}

\begin{equation}
\label{one_4}
{x \in \mathbb{Z}_{+}^{|\mathscr{R}|}},
\end{equation}

\begin{equation}
\label{one_5}
{y \in \mathbb{Z}_{+}^{|\mathscr{R}| \times |\mathscr{R}|}}.
\end{equation}

\par
The objective function (\ref{one_1}) is to maximize the overall profit which equals the difference of total revenue and total holding cost. The constraint in equation (\ref{one_2}) ensures that the number of total vehicles are not exceeded the capacity which can be easily estimated from historical data. The constraints in equation (\ref{one_3}) guarantee that each location must satisfy the customer demand. Constraints (\ref{one_4}) and (\ref{one_5}) are the types of decision variables.

\par
Although deterministic model is capable of tackling the optimization model in a simple way, the average demands for model may lead to optimal solution with high risk even infeasible. Additionally, it is worth noting that the objective function (\ref{one_1}) is a piece-wise linear function, therefore, it is required to reformulated to a linear function before solving.

\subsection{Two-Stage Stochastic Programming Model}
\par
The two-stage stochastic programming model of car-sharing problem can be formulated as follows.

\begin{equation}
\label{two_1}
\text{max} \sum_{s \in \mathscr{S}} P_{s}\left[\sum_{i \in \mathscr{R}} \min \left(x_{i}+\sum_{j \in \mathscr{R}} y_{i,j}^{s}, d_{i}^{s}\right) * r_{i} - \sum_{i \in \mathscr{R}}\left(h_{i} * x_{i}+\sum_{j \in \mathscr{R}} t_{i,j} * y_{i,j}^{s} \right) \right]
\end{equation}

\textbf{s.t.}

\begin{equation}
\label{two_2}
\sum_{i \in \mathscr{R}} x_{i} \leqslant \mathcal{C},
\end{equation}

\begin{equation}
\label{two_3}
\sum_{j \in \mathscr{R}} y_{i,j}^{s} \leq x_{i} \quad \forall i \in \mathscr{R}, \forall s \in \mathscr{S},
\end{equation}

\begin{equation}
\label{two_4}
{x \in \mathbb{Z}_{+}^{|\mathscr{R}|}},
\end{equation}

\begin{equation}
\label{two_5}
{y \in \mathbb{Z}_{+}^{|\mathscr{R}| \times |\mathscr{R}| \times |\mathscr{S}|}}.
\end{equation}

\par
Similar as one-stage SP model, the objective function in equation (\ref{two_1}) is to maximize the overall profit, which is denotes by the difference of revenue and overall cost (the summation of holding cost and moving/transferring cost). Constraint (\ref{two_2}) denotes the number of vehicle cannot exceed the capacity of car firm. Constraint (\ref{two_3}) implies that the sum number of cars that moving from location i to each location should not exceed the number of car at location i. Constraints (\ref{two_4}) and (\ref{two_5}) describe the type of decision variables.

\section{Approach}
\subsection{Model Reformulation}
Unlike the deterministic model which can be solved by off-the-shelf commercial solvers effectively. Normally, the two-stage SP model required reformulation since the continuous probability distribution contains infinite scenarios. In this paper, we utilize the sample average approximation (SAA)\cite{santoso2005stochastic} - a Monte Carlo method to reformulate the two-stage SP model. There are a variety of variant SAA approaches\cite{geyer1992constrained, mak1999monte,plambeck1996sample, shapiro1998simulation} with different names. In order to reduce the computation, we consider a simplified edition. The the procedure of SAA can be summarized as follows.

\begin{algorithm}
\caption{Sampling Average Approximation }\label{euclid}
\label{al}
\hspace*{\algorithmicindent} 
\textbf{Input: probability distribution $\mathcal{P}$, number of sample $M$, size $N$, two-stage SP model $z^{*}=\min _{x \in X} c^{T} x+\mathbb{E}_{\mathcal{P}}[Q(x, \xi(\omega))]$}  
\\
\hspace*{\algorithmicindent} 
\textbf{Output: the optimal value}

\begin{algorithmic}[1]

\State $k \gets 0$
\While{$k \leqslant M$}
\State $k \gets k+1$
\State a sample $\omega^{1}, \omega^{2},...,\omega^{n} $  of N scenario is generated according to $\mathcal{P}$;
\State reformulate the model as $z_{N}=\min _{x \in X} c^{T} x+\frac{1}{N} \sum_{n=1}^{N} Q\left(x, \xi\left(\omega^{n}\right)\right)$;
\State solve the model and get optimal value $z_{N}^{k}$ and optimal solution $\hat{x}^{k}$;

\EndWhile
\State\textbf{end while}

\State The optimal value is $\overline{z}_{N}=\frac{1}{M} \sum_{m=1}^{M} z_{N}^{m}$

\end{algorithmic}
\end{algorithm}

\par
Notice that the reformulation model in SAA, the objective function becomes

\begin{equation}
\text{max} \quad \mathbb{N}^{-1}\left[\sum_{i \in \mathscr{R}} \min \left(x_{i}+\sum_{j \in \mathscr{R}} y_{i,j}^{s}, d_{i}^{s}\right) * r_{i} - \sum_{i \in \mathscr{R}}\left(h_{i} * x_{i}+\sum_{j \in \mathscr{R}} t_{i,j} * y_{i,j}^{s} \right) \right]
\end{equation}

\par
where $\mathbb{N}$ is the number of scenarios. Additionally, the objective function is still a non-linear objective function. We introduce the auxiliary variable to transform the non-linear objective function to a linear type. Then the two-stage SP model becomes

\begin{equation}
\label{three_1}
\text{max} \quad \mathbb{N}^{-1}\left[\sum_{i \in \mathscr{R}} \sum_{s \in \mathscr{S}} r_{i} * f_{i}^{s} - \sum_{i \in \mathscr{R}}\left(h_{i} * x_{i}+\sum_{j \in \mathscr{R}} t_{i,j} * y_{i,j}^{s} \right) \right]
\end{equation}

\textbf{s.t.} 

\begin{equation}
\label{three_2}
f_{i}^{s} \leq x_{i} + \sum_{j \in \mathscr{R}} y_{i,j}^{s} \quad \forall i \in \mathscr{R}, \forall s \in \mathscr{S},
\end{equation}

\begin{equation}
\label{three_3}
f_{i}^{s} \leq d_{i}^{s} \quad \forall i \in \mathscr{R}, \forall s \in \mathscr{S},
\end{equation}

\begin{equation}
\label{three_4}
\sum_{i \in \mathscr{R}} x_{i} \leqslant Capacity,
\end{equation}

\begin{equation}
\label{three_5}
\sum_{j \in \mathscr{R}} y_{i,j}^{s} \leq x_{i} \quad \forall i \in \mathscr{R}, \forall s \in \mathscr{S},
\end{equation}

\begin{equation}
\label{three_6}
{x \in \mathbb{Z}_{+}^{|\mathscr{R}|}},
\end{equation}

\begin{equation}
\label{three_7}
{y \in \mathbb{Z}_{+}^{|\mathscr{R}| \times |\mathscr{R}|}}.
\end{equation}

\subsection{Model Decomposition}
\par 
After the final reformulation, the two-stage SP model becomes a very large-scale deterministic model, for example, if we consider 50 locations and 1000 scenarios, the number of second-stage decision variable will be 50*50*1000 = 2,500,000. To solve large-scale model effectively, decomposition algorithm is required. In this section, we introduce Benders decomposition\cite{benders2005partitioning} to solve the problem. For convenience, in the following, we neglect the constant N. Then we divide the reformulated model into master problem (MP)  
\begin{equation}
    \text{max} \sum_{i \in \mathscr{R}}\sum_{s \in \mathscr{S}} r_{i}* f_{i}^{s} - \sum_{i \in \mathscr{R}} h_{i} * x_{i} + \theta
\end{equation}
and subproblem (SP) in the dual form
\begin{equation}
    \text{min} \sum_{i \in \mathscr{R}}\sum_{s \in \mathscr{S}} (\bar f_{i}^{s} - \bar x_{i} - d_{i}^{s}) * u_{i} - \sum_{j \in \mathscr{R}} \bar x_{j} * v_{j}
\end{equation}

\textbf{s.t.} 

\begin{equation}
    u_{i} - v_{j} \leq t_{i,j} \quad\quad \forall i,j \in \mathscr{R}
\end{equation}

$u_{i}$ and $v_{j}$ are the dual variables of SP, while $\bar f_{i}^{s}$ and $\bar x_{i}$ are the fixed values that are determined by MP. During each iteration in MP, the values are adjusted and assigned to SP. Finally, the algorithm can be summarized as follows.

\begin{algorithm}
\caption{Benders Decomposition for Two-Stage SP Car-Sharing Model}
\label{al}
\hspace*{\algorithmicindent} 
\textbf{Input: $\mathcal{MP, SP}, \xi$}  \\
\hspace*{\algorithmicindent} 
\textbf{Output: the optimal solution }
\begin{algorithmic}[1]
\State $\mathcal{UB} \gets +\infty, \mathcal{LB} \gets -\infty$;
\While{$\mathcal{UB} - \mathcal{LB} \geq \xi$}
\State {given the fixed value $\bar f$ and $\bar x$ solve the $\mathcal{SP}$ model}
\If{$\mathcal{SP}$ is unbounded}
\State get ray($u^{*}$, $v^{*}$) and add cut $u^{*}*(f-x-d)$ - $v^{*}*x \leq 0$ to $\mathcal{MP}$ 
\ElsIf{$\mathcal{SP}$ is optimal}
\State get point($u^{*}$, $v^{*}$) and add cut $u^{*}*(f-x-d)$ - $v^{*}*x \leq \theta$ to $\mathcal{MP}$ 
\Else
\State the original model is infeasible
\EndIf
\State \textbf{end if}

\State solve the $\mathcal{MP}$ model
\State update $\mathcal{LB} \gets$  value of $\mathcal{MP}$
\EndWhile  \label{roy's loop}
\State \textbf{end while}
\State return either $\mathcal{LB}$ or $\mathcal{UB}$ as the optimal value
\end{algorithmic}
\end{algorithm}
\section{Numerical Experiment}
\textbf{Experiment Setup}. All the algorithms (KDE, SAA and BD) are implemented using Python 3.7, the mathematical models are solved by Gurobi 8.1 under the platform Intel i7, 16GB RAM, Windows 10. 

\par
\textbf{Experiment Design}. We devise a group of experiments. After the data preprocessing and distribution estimation by KDE, firstly, we validate the running times and expected profits based on different numbers of scenarios, additionally, we compare the outcomes deterministic model with two-stage SP model. Secondly, we compare the results yielded from non-parametric approach KDE with several parametric distributions like Gaussian distribution in terms of expected profit. The above experiments are based on training sets, finally, we compare the expected values obtained from training sets with testing sets. Specifically, we fix the first-stage decision variables by the outcomes that yields from training sets, then compute the overall expected values that the demands are from testing sets.

\subsection{Data Pre-possessing}
The data sets are from New York taxi trip\footnote{https://www1.nyc.gov/site/tlc/about/tlc-trip-record-data.page}, we collected three years (July 2016 - June 2019) green taxi trip records as the data source which is archived by month. Then we split the three years data sets into training set (from July 2016 to December 2018) and testing set (from January 2019 to June 2019), each data set involves thousands of naive one-trip records with a complex structure. Take the data set 2018-01 for example, it contains 793,529 records and 19 attributes. It is worth noting that the deterministic parameters in our SP model like $r_{i}$ (revenue) and $t_{ij}$ (transferring cost) can be estimated from the data set easily. For convenience, in the following experiments, the revenue per car is set to \$100, the transferring cost is set to \$5 by rough estimation, the number of available vehicles is set to 15,000, and the holding cost is assumed to follow the Gaussian distribution with parameters $\mathcal{N}(20, 9)$. Additionally, in this data set the whole New York city is divided into 259 different locations, we picked 20 locations with highest demands, which are aggregated by days (i.e. 914 days for training set and 181 days for testing set). Meanwhile, the pickup location and drop off location names are mapped as location ID stored in PULocationID and DOLocationID in the data set respectively. The New York city location division information details can be found via https://data.world/nyc-taxi-limo/taxi-zone-lookup.

\subsection{Probability Distribution and Sampling Results}
After using KDE, the demands probability distributions of each location is illustrated in the following figure, in which the bar plots denote the primitive demand while the curves are the approximate distributions for the locations derived by KDE.

\begin{figure}[H]
\centering
\includegraphics[scale=0.35]{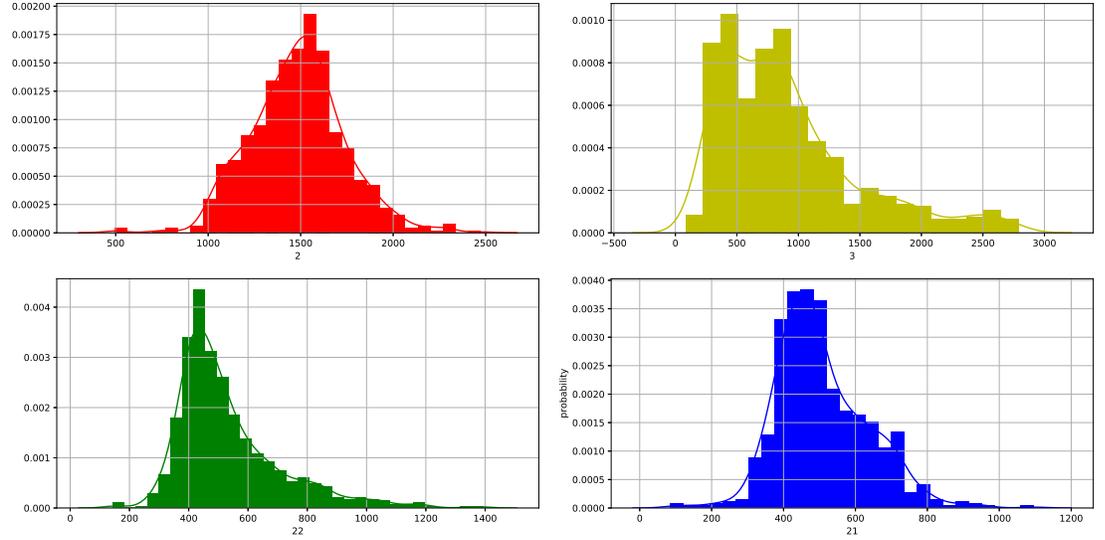}
\caption{Demand Distributions with Unimodal Type on Training Sets}
\label{fig:unimodal}
\end{figure}
Among the top 20 demands locations, there are mainly two types of distribution, one is unimodal type, which can be seen from Figure \ref{fig:unimodal}. The other type which represents the most locations is bimodal type. This can be seen from  the following figure. 

\par
In the first type, a specific functional form for the density model such as Gaussian distribution can be assumed, in other words, parametric methods can be applied on these scenarios. While in the second type, the particular form of parametric functions are unable to provide the appropriate representation of the real density. In such cases, we must consider using non-parametric approaches such as KDE.

\begin{figure}[H]
\centering
\includegraphics[scale=0.33]{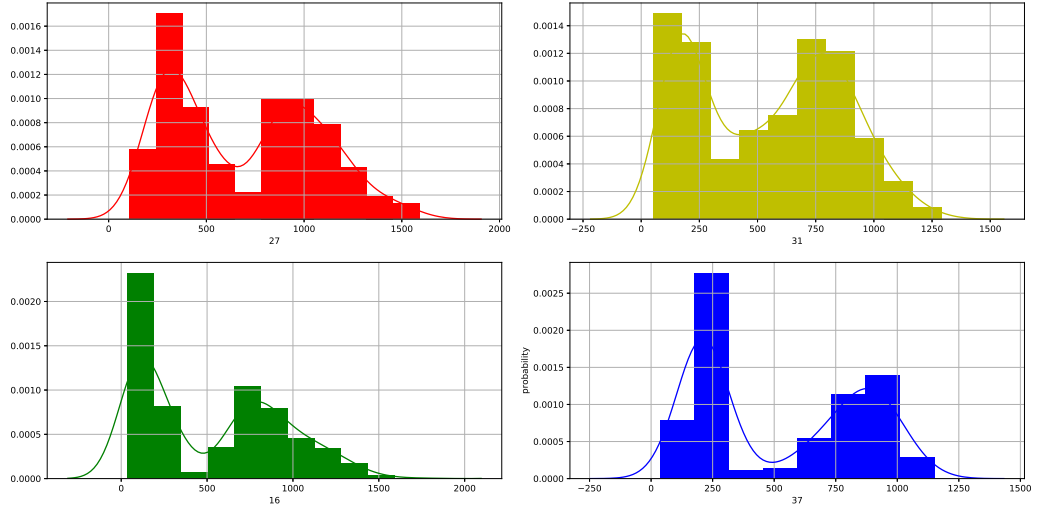}
\caption{Demand Distributions with Bimodal Type on Training Sets}
\label{fig:bimodal}
\end{figure}

\par
Most of the parametric methods may work well in the unimodal distributions, but cannot achieve the same goal for bimodal distributions. That is why KDE approach is introduced in this work. 

\subsection{Stochastic Model vs. Deterministic Model Results}
In this experiment, We generate 5 groups of scenarios for SP model based on the probability distributions that are derived from KDE. The numbers of scenarios are 20, 50, 100, 200 and 500, each group runs 10 times. Additionally, we consider deterministic model using the average demands that are calculated from training set (average demand of 919 days) and testing set (average demand of 181 days). The average objective value and time elapse can be seen in the table below.

\begin{table}[h]
\centering
\begin{tabular}{ c c c }
\hline
Number of Scenario & Objective Value & Time Elapse (s) \\ \hline
20       & \$1,477,845                               & 2.73            \\ 
50       & \$1,487,606                               & 6.87            \\ 
100      & \$1,475,688                               & 10.89           \\ 
200      & \$1,484,367                               & 21.73           \\ 
500      & \$1,469,642                               & 53.12           \\ \hline
deterministic (average on training set)     & \$1,325,723                                & 0.24          \\ \
deterministic (average on testing set)     & \$1,017,054                               & 0.24          \\ \hline
\end{tabular}
\caption{Average Objective Value and Time Elapse under Different Number of Scenarios}
\label{table: avg}
\end{table}

\par Based on the experimental results, we come to conclude that the two-stage SP model is able to yield more outcomes than the deterministic model. While by average demands, the overall profit on the training set is more that the one on the testing set. Meanwhile, as we discussed in the beginning of this section, as the number of scenario increases, the time elapses grows by approximate linear increment.

\subsection{Validations on Parametric Distributions}

\par 
Besides the non-parametric approach, we also use several popular parametric distributions (such as Gaussian, lognormal, Laplace and Exponential distributions) in terms of average expected value. Among these distributions, we found that the SP model based on the exponential and lognormal distributions are infeasible. The reason is the high average demand in a specific location will yield the sampling that with extreme high demand which lead the SP model infeasible.

\begin{table}[]
\centering
\begin{tabular}{c c c}
\hline
Number of Scenario & Outcome(Laplace)     & Outcome(Gaussian)    \\ \hline
20       & \$1,467,117 & \$1,425,569 \\ 
50       & \$1,422,868 & \$1,402,279 \\ 
100      & \$1,417,811 & \$1,417,403 \\ 
200      & \$1,406,112 & \$1,412,343 \\ 
500      & \$1,406,103 & \$1,398,546 \\ \hline
\end{tabular}
\caption{Average Objective Value under Different Probability Distributions}
\end{table}

\par 

\par
As can be seen from the table, the overall profit yielded from Laplace distribution is slightly better than the one yielded from Gaussian distribution. However both of the  parametric approaches are inferior to the one from KDE in terms of the overall profit.

\subsection{Solutions Comparisons on Testing Sets}

In the two-stage SP model, solutions involves two parts, first-stage decision variables which denotes the numbers of cars that are placed at each location before demands realize. Second-stage decision variables which denote the number of cars that are moving between locations. In this experiment, we use the values of first-stage decision variables that are derived from two-stage SP model that are constructed based on training model, to validate the overall profit of two-stage SP model that are constructed based on testing sets (6 months, 181 days).

\begin{figure}[htbp]
\centering
\includegraphics[scale=0.8]{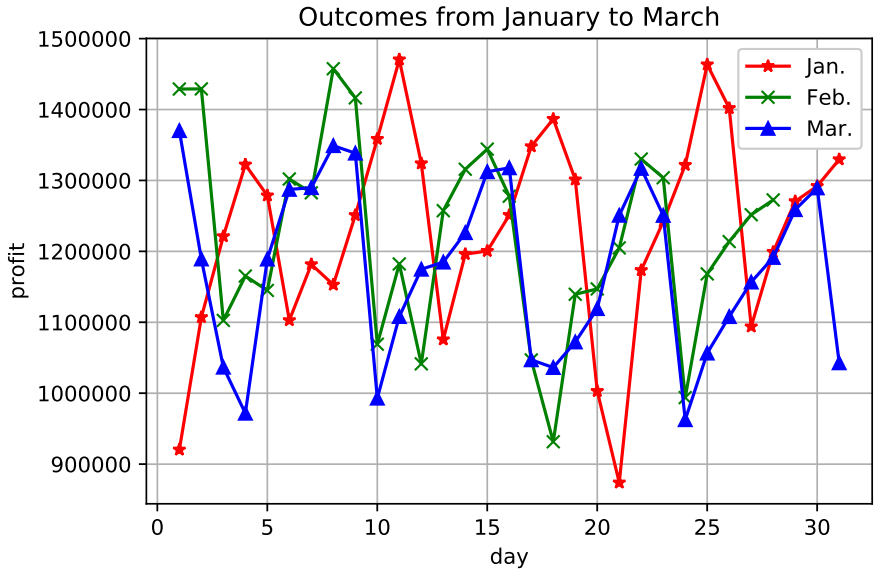}
\caption{Overall Profit on Testing Sets 2019-01 to 2019-03}
\end{figure}

\begin{figure}[htbp]
\centering
\includegraphics[scale=0.8]{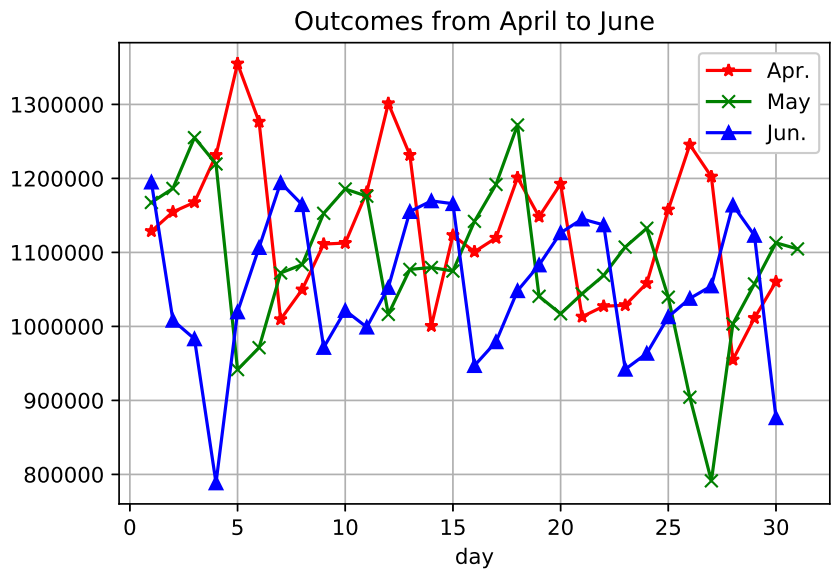}
\caption{Overall Profit on Testing Sets 2019-04 to 2019-06}
\end{figure}

\par
As can be seen from the results, using the values first-stage decision variables, the outcomes gap on training sets and testing sets is from 0.35\% to 46.57\%, the average gap is 22.1\%. Meanwhile, compared to Table \ref{table: avg}, we come to conclude that most profits that yielded from SP model on the testing sets are higher that the profit that yielded from deterministic model on the testing sets.

\section{Conclusions and Future Work}
\par
In this paper, we propose a framework that involves kernel density estimation to predict the location demands, and a two-stage stochastic programming model to solve the car-sharing problem under demand uncertainty. In more real world, the demand distribution would be time variant and evolves gradually (or the parameters of distribution vary at least), which renders the data-driven system outdated and leads to deteriorates the resulting solution quality\cite{ning2019optimization}. In order to describe this evolution in a more precise way, we will investigate Bayesian learning which focus on posterior probability distribution that is based on prior probability distribution and the likelihood of current data. Namely, we will explore the dynamic data-driven stochastic programming model for car-sharing problem.

\par
Additionally, in our work, the proposed framework treats the location demands by days. For some real-time applications, however, the daily demand should be considered as time-series data, which would be handled by time-series prediction machine learning algorithm. We will explore this topic in our future works.

\par
Meanwhile, in this paper, for convenience, some other factors we do not consider. For example, we do not consider the capacity of locations, and the route condition of balancing which may lead different transportation costs. Later on, we will extend the two-stage SP model to a more practical one.

\bibliographystyle{model1-num-names}
\bibliography{elsarticle-template-1-num.bbl}





\end{document}